\documentstyle[amssymb,amsfonts]{amsart}

\def\R{{\hbox{\bf R}}}
\def\C{{\hbox{\bf C}}}

 at 10 true pt

\def\T{{\hbox{\bf T}}}

\def\allt#1{%
\smash{
\vtop{%
     \ialign{%
        ##\crcr
        $\hfil\displaystyle{\tilde \forall}\hfil$\crcr%
        \noalign{\kern1.5pt\nointerlineskip}
        $\hfil\!\!#1\hfil$\crcr\noalign{\kern1.5pt}
        }
       }
      } \hbox{$\vphantom{#1}$}
     }

\def\be#1{\begin{equation}\label{#1}}
\def\bas{\begin{align*}}
\def\eas{\end{align*}}
\def\bi{\begin{itemize}}
\def\ei{\end{itemize}}

\def\eps{\varepsilon}

\def \endprf{\hfill  {\vrule height6pt width6pt depth0pt}\medskip}
\def\emph#1{{\it #1}}
\def\textbf#1{{\bf #1}}

\theoremstyle{plain}
  \newtheorem{theorem}[subsection]{Theorem}

\theoremstyle{remark}

\theoremstyle{definition}

\include{psfig}

\begin{document}

\title[Recent progress on Restriction]{Some recent progress on the Restriction conjecture}

\author{Terence Tao}
\address{Department of Mathematics, UCLA, Los Angeles CA 90095-1555}
\email{tao@@math.ucla.edu}


\begin{abstract} 
We survey recent developments on the Restriction conjecture.
\end{abstract}

\maketitle

\section{Introduction}

The purpose of this paper is describe the state of progress on the restriction problem in harmonic analysis, with an emphasis on the developments of the past decade or so on the Euclidean space version of these problems for spheres and other hypersurfaces.  As the field is quite large and has so many applications, it will be impossible to completely survey the field, but we will try to at least give the main ideas and developments in this area.

The restriction problem are connected to many other conjectures, notably the Kakeya and Bochner-Riesz conjectures, as well as PDE conjectures such as the local smoothing conjecture.  For reasons of space, we will not be able to discuss all these connections in detail, but refer the reader to \cite{wolff:survey}, \cite{Bo}, \cite{tao:elesc}, \cite{tao:notices}; for the connection between restriction and Bochner-Riesz see \cite{feff:thesis}, \cite{carbery:parabola}, \cite{borg:stein}, \cite{borg:kakeya}, \cite{tao:boch-rest}, \cite{lee:2}, etc.

\section{The restriction problem}\label{sec:restriction}

Fix $n \geq 2$; in our discussion all the constants $C$ are allowed to depend\footnote{The question on how quantifying the exact dependence of the constants here on the dimension $n$ as $n \to \infty$ is however an interesting problem, although to my knowledge there are not many results in this direction at present.} on $n$.

We begin by discussing the restriction problem.  Historically, this problem originated by studying the Fourier transform of $L^p$ functions in Euclidean space $\R^n$ for some $n \geq 1$, although it was later realized that this problem also arises naturally in other contexts, such as non-linear PDE and in the study of eigenfunctions of the Laplacian.  

If $f$ is an $L^1(\R^n)$ function, then the Riemann-Lebesgue lemma implies that the Fourier transform $\hat f$, defined by
$$ \hat f(\xi) := \int_{\R^n} e^{-2\pi i x \cdot \xi} f(x)\ dx$$
is a continuous bounded function on $\R^n$ which vanishes at infinity.  In particular, we can meaningfully restrict this function to any subset $S$ of $\R^n$, creating a continuous bounded function $\hat f|_S$ on $S$.

On the other hand, if $f$ is an arbitrary $L^2(\R^n)$ function, then the Fourier transform $\hat f$ can be any function in $L^2(\R^n)$, and in particular there is no meaningful way to restrict it to any set $S$ of zero measure.  

Between these two extremes, one may ask what happens to the Fourier transform of a function $f$ in $L^p(\R^n)$, where $1 < p < 2$.  Certainly we do not expect the Fourier transform $\hat f$ to be continuous or bounded, and it is easy to construct examples of $L^p$ functions which have an infinite Fourier transform at one point.  In fact, it is easy to create such a function which is infinite on an entire hyperplane; for instance, the function
$$ f(x) := \frac{1}{1 + |x_1|}$$
where $x_1$ is the first co-ordinate of $x$, lies in $L^p$ for every $p>1$, but has an infinite Fourier transform on every point on the hyperplane $\{ \xi \in \R^n: \xi_1 = 0 \}$.

On the other hand, from the Hausdorff-Young inequality we see that $\hat f$ lies in the Lebesgue space $L^{p'}(\R^n)$, where $1/p + 1/p' = 1$.  Thus $\hat f$ can be meaningfully restricted to every set $S$ of positive measure.  

This leaves open the question of what happens to sets $S$ which have zero measure but which are not contained in hyperplanes.  In 1967 Stein made the surprising discovery that when such sets contain sufficient ``curvature'', that one can indeed restrict the Fourier transform of $L^p$ functions for certain $p > 1$.  This lead to the \emph{restriction problem} \cite{stein:problem}: for which sets $S \subseteq \R^n$ and which $1 \leq p \leq 2$ can the Fourier transform of an $L^p(\R^n)$ function be meaningfully restricted?

There are of course infinitely many such sets to consider, but we shall focus our attention here on sets $S$ which are hypersurfaces\footnote{For surfaces of lower dimension, see \cite{christ:thesis}, \cite{prestini}, \cite{mock:habil}, \cite{banner}; for fractal sets in $\R$, see \cite{mock:fractal}, \cite{sj2}, \cite{mock:habil}; for surfaces in finite field geometries, see \cite{mock:tao}; for the restriction theory of the prime numbers, see \cite{green}.}, or compact subsets of hypersurfaces.  In particular, we shall be interested\footnote{It is easy to see, using the symmetries of the Fourier transform, that the restriction problem for a set $S$ is unaffected by applying any translations or invertible linear transformations to the set $S$, so we can place the sphere, paraboloid, and cone in their standard forms \eqref{eq:sphere}, \eqref{eq:paraboloid}, \eqref{eq:cone} without loss of generality.} in the \emph{sphere}
\begin{equation}\label{eq:sphere}
S_{sphere} := \{ \xi \in \R^n: |\xi| = 1\},
\end{equation}
the \emph{paraboloid}
\begin{equation}\label{eq:paraboloid}
S_{parab} := \{ \xi \in \R^n: \xi_n = \frac{1}{2} |\underline \xi|^2 \},
\end{equation}
and the \emph{cone}
\begin{equation}\label{eq:cone}
S_{cone} := \{ \xi \in \R^n: \xi_n = |\underline \xi| \},
\end{equation}
where $\xi = (\underline{\xi}, \xi_n) \in \R^{n-1} \times \R \equiv \R^n$, and we always take $n \geq 2$ to avoid trivial situations.
These three surfaces are model examples of hypersurfaces with curvature\footnote{One could also consider cylinders such as $S^{k-1} \times \R^{n-k} \subset \R^n$, but it turns out that the restriction theory for these surfaces is identical to that of the sphere $S^{k-1}$ inside $\R^k$.}, though of course the cone differs from the sphere and paraboloid in that it has one vanishing principal curvature.  These three hypersurfaces also enjoy a large group of symmetries (the orthogonal group, the parabolic scaling and Gallilean groups, and the Poincare group, respectively).  Also, these three hypersurfaces are related via the Fourier transform to solutions to certain familiar partial differential equations, namely the Helmholtz equation, Schr\"odinger equation, and wave equation; however we will not focus on applications to PDE in this paper.

\section{Restriction estimates and extension estimates }\label{sec:rest-extend}

Let $S$ be a compact subset (but with non-empty interior) of one of the above surfaces $S_{sphere}$, $S_{parab}$, $S_{cone}$.   We endow $S$ with a canonical measure $d\sigma$ - for the sphere, this is surface measure, for the parabola, it is the pullback of the $n-1$-dimensional Lebesgue measure $d\underline{\xi}$ under the projection map $\xi \mapsto \underline \xi$, while for the cone the pullback of $d\underline{\xi}/|\xi|$ is the most natural measure (as it is Lorentz-invariant.  In order to restrict the Fourier transform of an $L^p(\R^n)$ function to $S$, it will suffice to prove an \emph{a priori} ``restriction estimate'' of the form
\begin{equation}\label{eq:rpq}
 \| \hat f|_S \|_{L^q(S; d\sigma)} \leq C_{p,q,S} \| f \|_{L^p(\R^n)}
\end{equation}
for all Schwartz functions $f$ and some $1 \leq q \leq \infty$, since one can then use density arguments to obtain a continuous restriction operator from $L^p(\R^n)$ to $L^q(S; d\sigma)$ which extends the map $f \mapsto \hat f|_S$ for Schwartz functions.  When the set $S$ has sufficient symmetry (e.g. if $S$ is the sphere), this implication can in fact be reversed, using Stein's maximal principle \cite{stein:maximal}; if there is no bound of the form\footnote{Indeed, it suffices for the weak-type estimate from $L^p(\R^n)$ to $L^{p,\infty}(S; d\sigma)$ to fail.  See \cite{stein:maximal}; similar ideas arise in the factorization theory of Nikishin and Pisier.} \eqref{eq:rpq}, then one can construct functions $f \in L^p(\R^n)$ whose Fourier transform is infinite almost everywhere in $S$.

We will tend to think of $\R^n$ as representing ``physical space'', whose elements will be denoted names such as $x$ and $y$, while $S$ lives in ``frequency space'', and whose elements will be denoted names such as $\xi$ or $\omega$.  For the PDE applications it is sometimes convenient to think of $\R^n$ as a spacetime $\R^{n-1} \times \R := \{ (x,t): x \in \R^{n-1}, t \in \R \}$ (with the frequency space thus becoming spacetime frequency space $\{ (\xi,\tau): \xi \in \R^{n-1}, \tau \in \R \}$), but we will avoid doing so here.

It is thus of interest to see for which sets $S$ and which exponents $p$ and $q$ one has estimates of the form \eqref{eq:rpq}; henceforth we assume our functions $f$ to be Schwartz.  We denote by $R_S(p \to q)$ the statement that \eqref{eq:rpq} holds for all $f$.  From our previous remarks we thus see that $R_S(1 \to q)$ holds for all $1 \leq q \leq \infty$, while $R_S(2 \to q)$ fails for all $1 \leq q \leq \infty$; the interesting question is then what happens for intermediate values of $p$.  If $S$ is compact, then an estimate of the form $R_S(p \to q)$ will automatically imply an estimate $R_S(\tilde p \to \tilde q)$ for all $\tilde p \leq p$ and $\tilde q \leq q$ by the Sobolev and H\"older inequalities.  Thus the aim is to increase the size of $p$ and $q$ for which $R_S(p \to q)$ holds by as much as possible.

A simple duality argument shows that the estimate \eqref{eq:rpq} is equivalent to the ``extension estimate''
\begin{equation}\label{eq:rpq-dual}
 \| (F d\sigma)^\vee \|_{L^{p'}(\R^n)} \leq C_{p,q,S} \| F \|_{L^{q'}(S; d\sigma)}
\end{equation}
for all smooth functions $F$ on $S$, where $(F d\sigma)^\vee$ is the inverse Fourier transform of the measure $Fd\sigma$:
$$ (Fd\sigma)^\vee(x) := \int_S F(\xi) e^{2\pi i \xi \cdot x} d\sigma(\xi).$$
Indeed, the equivalence of \eqref{eq:rpq} and \eqref{eq:rpq-dual} follows from Parseval's identity 

$$\int_{\R^n} (F d\sigma)^\vee(x) \overline{f(x)}\ dx = \int_S F(\xi) \hat{\overline{f(\xi)}}\ d\sigma(\xi)$$
and duality.  If we use $R^*_S(q' \to p')$ to denote the statement that the estimate \eqref{eq:rpq-dual} holds, then $R^*_S(q' \to p')$ is thus equivalent to $R_S(p \to q)$.

Note that because $F$ is smooth, it is possible to use the principle of stationary phase (see e.g. \cite{stein:large}) to obtain asymptotics for $(F d\sigma)^\vee$.  However, such asymptotics depend very much on the smooth norms of $F$, not just on the $L^{q'}(S)$ norm, and so do not imply estimates of the form \eqref{eq:rpq-dual} (although they can be used to provide counterexamples).  Thus one can think of extension estimates as a more general way than stationary phase to control oscillatory integrals, applicable in situations where the amplitude function $F(\xi)$ has magnitude bounds but no smoothness properties.

The extension formulation \eqref{eq:rpq-dual} also highlights the connection between this problem and partial differential equations.  For instance, consider a solution $u(t,x): \R \times \R^n \to \C$ to the free Schr\"odinger equation
$$ i\partial_t u + \Delta u = 0$$
with initial data $u(0,x) = u_0(x)$.  This has the explicit solution
$$ u(t,x) = \int e^{2\pi i (x \cdot \xi + 2\pi t |\xi|^2)} \hat u_0(\xi)\ d\xi,$$
or equivalently
$$ u = (F d\sigma)^\vee$$
where $d\sigma := d\xi \delta(\tau - 2\pi |\xi|^2)$ is (weighted) surface measure on the paraboloid $\{ (\tau,\xi) \in \R \times \R^n: \tau = 2 \pi |\xi|^2 \}$, and $F$ is the function $\hat u_0(\xi)$ restricted to the paraboloid.  Thus, estimates of the form $\R^*_S(q' \to p')$ when $S$ is the paraboloid in $\R \times \R^n$ to control certain spacetime norms of solutions to the free Schr\"odinger equation.  Somewhat similar connections exist between the cone \eqref{eq:cone} (in $\R \times \R^n$) and solutions to the wave equation $u_{tt} - \Delta u = 0$, or between the sphere \eqref{eq:sphere} and solutions to the Helmholtz equation $\Delta u + u = 0$.  We will not pursue these connections further here, but see for instance \cite{strichartz:restrictionquadratic} and the (numerous) papers descended from that paper.  (Some other connections between restriction estimates and PDE-type estimates are summarized in \cite{tv:cone2} and the references therein; for the Helmholtz equation, see for instance \cite{ruiz}.)

\section{Necessary conditions}\label{sec:necessary}

We will use the extension formulation \eqref{eq:rpq-dual} to develop some necessary conditions in order for $\R^*_S(q' \to p')$ to hold.  First of all, by setting $F \equiv 1$ we clearly see that we must have $(d\sigma)^\vee \in L^{p'}(\R^n)$ as a necessary condition.  In the case of the sphere \eqref{eq:sphere}, the Fourier transform $(d\sigma)^\vee(x)$ decays in magnitude like $(1 + |x|)^{-(n-1)/2}$ (as can be seen either by stationary phase, or by the asymptotics of Bessel functions), and so we obtain the necessary condition\footnote{There does not seem to be any hope for any weak-type endpoint estimate at $p' = 2n/(n-1)$, see \cite{bcss}.} $p' > 2n/(n-1)$, or equivalently $p < 2n/(n+1)$.  A similar computation gives the same constraint $p' > 2n/(n-1)$ for the paraboloid \eqref{eq:paraboloid}, while for the cone the asymptotics are slightly different, giving the condition $p' > 2(n-1)/(n-2)$.

Let $F$ be a smooth function on $S$ with an $L^\infty$ norm of at most 1.  Since $Fd\sigma$ is pointwise dominated by $d\sigma$, it seems intuitive that $(Fd\sigma)^\vee$ should be ``smaller'' than $(d\sigma)^\vee$.  Thus one should expect the above necessary conditions to in fact be sufficient to obtain the estimate $R^*_S(\infty \to p')$.  For completely general sets $S$, this assertion is essentially the \emph{Hardy-Littlewood majorant conjecture}; it is true when $p'$ is an even integer by direct calculation using Plancherel's theorem, but is false for other values of $p'$ (a ``logarithmic'' failure was established by Bachelis in the 1970s; a more recent ``polynomial'' failure has been established independently by Mockenhaupt and Schlag (private communication) and Green and Ruzsa (private communication).  See \cite{mock:habil} for further discussion).  However, it may still be that the majorant conjecture is still true for ``non-pathological'' sets $S$ such as the sphere, paraboloid, and cone.

Another necessary condition comes from the \emph{Knapp example} \cite{tomas:restrict}, \cite{strichartz:restrictionquadratic}.  In the case of the sphere or paraboloid, we sketch the example as follows.  Let $R \gg 1$.  Then, by a Taylor expansion of the surface $S$ around any interior point $\xi_0$, we see that the surface $S$ contains a ``cap'' $\kappa \subset S$ centered at $\xi_0$ of diameter\footnote{We use $X \lesssim Y$ or $X = O(Y)$ to denote an estimate of the form $X \leq CY$ where $C$ depends on $S$, $p$, $q$, but not on functions such as $f$, $F$, or on parameters such as $R$.  We use $X \sim Y$ to denote the estimate $X \lesssim Y \lesssim X$.} $\sim 1/R$ and surface measure $\sim R^{-(n-1)}$ which is contained inside a disk $D$ of radius $\sim 1/R$ and thickness $\sim 1/R^2$, oriented perpendicular to the unit normal of $S$ at $x_0$.  Let $F$ be the characteristic function of this cap $\kappa$ (one can smooth $F$ out if desired, but this does not affect the final necessary condition), and let $T$ be the dual tube to the disk $D$, i.e. a tube centered at the origin of length $\sim R^2$ and thickness $\sim R$ oriented in the direction of the unit normal to $S$ at $x_0$.  Then $(F d\sigma)^\vee$ has magnitude $\sim \sigma(K) \sim R^{-(n-1)}$ on a large portion of $T$ (this is basically because for a large portion of points $x$ in $T$, the phase function $e^{2\pi i x \cdot \xi}$ is essentially constant on $K$).  In particular, we have
$$ \| (F d\sigma)^\vee \|_{L^{p'}(\R^n)} \gtrsim |T|^{1/p'} R^{-(n-1)} \sim R^{-(n-1)} R^{(n+1)/p'},$$
while we have
$$ \| F \|_{L^{q'}(S; d\sigma)} \lesssim |\kappa|^{-1/q'} \lesssim R^{-(n-1)/q'}.$$
Letting $R \to \infty$, we thus see that we need the necessary condition
$$ \frac{n+1}{p'} \leq \frac{n-1}{q}$$
in order for $R^*_S(q' \to p')$ to hold.  (In the case of compact subsets of the paraboloid with non-empty interior, one can obtain the same necessary condition using the parabolic scaling $(\underline \xi, \xi_n) \mapsto (\lambda \underline \xi, \lambda^2 \xi_n)$.  For the full (non-compact) paraboloid, one can improve this to $\frac{n+1}{p'} = \frac{n-1}{q}$).  In the case of the cone, we can lengthen the cap $\kappa$ in the null direction (so that it now has measure $\sim R^{-(n-2)}$ and lives in a ``plate'' of length $\sim 1$, width $\sim 1/R$ and thickness $1/R^2$), which eventually leads to the stronger necessary condition $\frac{n}{p'} \leq \frac{n-2}{q}$; as before, this can be strengthened to $\frac{n}{p'} = \frac{n-2}{q}$ if one is considering the full cone \eqref{eq:cone} and not just compact subsets of it with non-empty interior.

One can formulate a Knapp counterexample for any smooth hypersurface; the necessary conditions obtained this way become stronger as the surface becomes flatter, and in the extreme case where the surface is infinitely flat (e.g. when it is a hyperplane), there are no estimates.

The \emph{restriction conjecture} for the sphere, paraboloid, and cone then asserts that the above necessary conditions are in fact sufficient.  In other words, for compact subsets of the sphere and paraboloid the conjecture asserts that $R^*_S(q' \to p')$ holds when $p' > 2n/(n-1)$ and $\frac{n+1}{p'} \leq \frac{n-1}{q}$, while for compact subsets of the cone the conditions become $p' > 2(n-1)/(n-2)$ and $\frac{n}{p'} \leq \frac{n-2}{q}$ (i.e. they match the numerology of the sphere and paraboloid in one lower dimension).  This conjecture has been solved for the paraboloid and sphere in two dimensions, and for the cone in up to four dimensions; see Figures 1 and 2 for a more detailed summary of progress on this problem.  The restriction problems for the three surfaces are related; the sharp restriction conjecture for the sphere would imply the sharp restriction estimate for the paraboloid, because one can parabolically rescale the sphere to approach the paraboloid; see \cite{tao:boch-rest}.  Also, using the method of descent, one can link the restriction conjecture for the cone in $\R^{n+1}$ with the restriction conjecture for the sphere or paraboloid in $\R^n$, although the connection here is not as tight (see \cite{tao:cone} for some further discussion).

\begin{figure}\label{fig:parab}
\begin{tabular}{|l|l|l|} \hline
Dimension & Range of $p$ and $q$ & \\
\hline
$n = 2$ & $q' = 2, p' \geq 8$ & Stein, 1967 \\
& $q' \geq (p'/3)'; p' > 4$ & Zygmund, 1974 \cite{zygmund} (best possible)\\
\hline
$n = 3$  & $q' = 2, p' \geq 6$ & Stein 1967 \\
  & $q' > (p'/2)', p' > 4$ & Tomas 1975 \cite{tomas:restrict} \\
  & $q' \geq (p'/2)', p' \geq 4$ & Stein 1975; Sj\"olin $\sim$ 1975 \\
  & $q', p' > 4 - \frac{2}{15}$ & Bourgain 1991 \cite{borg:kakeya} \\
  & $q', p' > 4 - \frac{2}{11}$ & Wolff 1995 \cite{wolff:kakeya} \\
  & $q' > 7/3; p' > 4 - \frac{2}{11}$ & Moyua, Vargas, Vega 1996 \cite{vargas:restrict} \\
  & $q' \geq (p'/2)'; p' > 4 - \frac{5}{27}$ & Tao, Vargas, Vega 1998 \cite{tvv:bilinear}\\

  & $q' > 170/77; p' > 4 - \frac{2}{9}$ & Tao, Vargas, Vega 1998 \cite{tvv:bilinear}\\
  & $q' \geq (p'/2)'; p' > 4 - \frac{8}{31}$ & Tao, Vargas 2000 \cite{tv:cone1}\\
  & $q' > 26/11; p' > 4 - \frac{2}{7}$ & Tao, Vargas 2000 \cite{tv:cone1}\\
  & $q' \geq (p'/2)'; p' > 4 - \frac{2}{3}$ & Tao 2003 \cite{tao:parabola}\\
  & $q' \geq (p'/2)'; p' > 3$ & (conjectured)\\
\hline
$n > 3$ & $q' > ((n-1)p'/(n+1))'; p' > \frac{2(n+1)}{n-1}$ & Tomas 1975 \cite{tomas:restrict} \\
& $q' \geq ((n-1)p'/(n+1))'; p' \geq \frac{2(n+1)}{n-1}$ & Stein 1975 \\
& $q',p' > \frac{2(n+1)}{n-1} - \eps_n$ & Bourgain 1991 \cite{borg:kakeya} \\
& $q',p' > \frac{2n^2 + n + 6}{n^2 + n - 1}$  & Wolff 1995 \cite{wolff:kakeya} \\
& $q' > \frac{2n^2 + n + 6}{n^2 + 3n + 1}; p' > \frac{2n^2 + n + 6}{n^2 + n - 1}$  & Moyua, Vargas, Vega 1996 \cite{vargas:restrict} \\
  & $q' \geq ((n-1)p'/(n+1))'; p' > \frac{2(n+2)}{n}$ & Tao 2003 \cite{tao:parabola}\\
& $q' \geq ((n-1)p'/(n+1))'; p' > \frac{2n}{n-1}$ & (conjectured)\\
\hline
\end{tabular}
\caption{Known results on the restriction problem $R_S(p \to q)$ (or $R^*_S(q' \to p')$) for the sphere and for compact subsets of the paraboloid.  (For the whole paraboloid, restrict the above exponents to the range $q' = (\frac{(n-1)p'}{n+1})'$).}
\end{figure}

\begin{figure}\label{fig:cone}
\begin{tabular}{|l|l|l|} \hline
Dimension & Range of $p$ and $q$ & \\
\hline
$n = 3$  & $q' \geq (p'/3)', p' \geq 6$ & Strichartz 1977 \cite{strichartz:restrictionquadratic} \\
& $q' \geq (p'/3)'; p' > 4$ & Barcelo, 1985 \cite{barcelo} (best possible)\\
\hline
$n = 4$ &  $q' \geq (p'/2)', p' \geq 4$ & Strichartz 1977 
\cite{strichartz:restrictionquadratic} \\
& $q' \geq (p'/2)'; p' > 3$ & Wolff, 2000 \cite{wolff:cone} (best possible)\\
\hline
$n > 4$ &  $q' \geq ((n-2)p'/n)', p' \geq \frac{2n}{n-2}$ & Strichartz 1977 
\cite{strichartz:restrictionquadratic} \\
& $q' \geq ((n-2)p'/n)'; p' > \frac{2(n+2)}{n}$ & Wolff, 2000 \cite{wolff:cone}\\
& $q' \geq ((n-2)p'/n)'; o' > \frac{2(n-1)}{n-2}$ & (conjectured)\\
\hline
\end{tabular}
\caption{Known results on the restriction problem $R_S(p \to q)$ (or $R^*_S(q' \to p')$) for compact subsets of the cone.  (For the whole cone, restrict the above exponents to the range $q' = (\frac{(n-2)p'}{n})'$).}
\end{figure}

\section{Local restriction estimates}\label{sec:local-rest}

We now begin discussing some of the tools used to prove the above restriction theorems.  The first key idea is to reduce the study of \emph{global} restriction theorems (where the physical space variable is allowed to range over all of $\R^n$), to that of \emph{local} restriction theorems (where the physical space variable is constrained to lie in a ball).

More precisely, for any exponents $p,q$, and any $\alpha \geq 0$, let $R_S(p \to q; \alpha)$ denote the statement that the localized restriction estimate
\begin{equation}\label{eq:local-rest}
  \| \hat f|_S \|_{L^q(S; d\sigma)} \leq C_{p,q,S,\alpha} R^\alpha \| f \|_{L^p(B(x_0,R))}
\end{equation}
holds for any radius $R \geq 1$, any ball $B(x_0,R) := \{ x \in \R^n: |x-x_0| \leq R\}$ of radius $R$, and any test function $f$ supported in $B(x_0,R)$.  Note that the center $x_0$ of the ball is irrelevant since one can translate $f$ by an arbitrary amount without affecting the \emph{magnitude} of $\hat f$.

Observe that estimates for higher $\alpha$ immediately imply estimates for lower $\alpha$ (keeping $p$, $q$, $S$ fixed).  Also, the local estimate $R_S(p \to q; 0)$ is clearly equivalent to the global estimate $R_S(p \to q)$ by a sending $R \to \infty$ and applying a limiting argument.  Finally, it is easy to prove estimates of this type for very large $\alpha$; for instance, for smooth compact hypersurfaces $S$ one has the estimate $R_S(p \to q; n/p')$ just from the H\"older inequality
$$ |\hat f(\xi)| \leq \| f\|_1 \leq C_p R^{n/p'} \| f \|_{L^p(B(x_0,R)}.$$
Thus the aim is to lower the value of $\alpha$ from the trivial value of $\alpha = n/p'$, toward the ultimate aim of $\alpha = 0$, at least when $p$ and $q$ lie inside the conjectured range of the Restriction conjecture.  (For other $p$ and $q$, the canonical counterexamples will give some non-zero lower bound on $\alpha$).

By duality, the local restriction estimate $R_S(p \to q; \alpha)$ is equivalent to the local extension estimate $R^*_S(q' \to p'; \alpha)$, which asserts that
\begin{equation}\label{eq:local-extension}
\| (F d\sigma)^\vee \|_{L^{p'}(B(x_0,R))} \leq C_{p,q,S,\alpha} R^\alpha \| F \|_{L^{q'}(S; d\sigma)}
\end{equation}
for all smooth functions $F$ on $S$, all $R \geq 1$, and all balls $B(x_0,R)$.

The uncertainty principle suggests that since the spatial variable has now been localized to scale $R$, the frequency variable can be safely blurred to scale $1/R$.  In the case where $S$ is a smooth compact hypersurface, this is indeed correct; the estimate \eqref{eq:local-rest} is equivalent to the estimate

\begin{equation}\label{eq:local-rest-blur}
  \| \hat f \|_{L^q(N_{1/R}(S))} \leq C_{p,q,S,\alpha} R^{\alpha-1/q} \| f \|_{L^p(B(x_0,R))}
\end{equation}
holding for all test functions\footnote{One can in fact remove the hypothesis that $f$ is supported on a ball, and replace \eqref{eq:local-rest-blur} with the corresponding global estimate where $f$ ranges over $\R^n$, provided that we have the Hausdorff-Young condition $q \leq p'$.  This is another manifestation of the uncertainty principle; it can be proven by using a smooth partition of unity to divide a global $f \in L^p(\R^n)$ into functions in $L^p(B(x_0,R))$ for various balls $B(x_0,R)$ and applying \eqref{eq:local-rest-blur} to each piece.  To sum, one can subdivide $N_{1/R}(S)$ into cubes of size $1/R$ and apply a local form of the Hausdorff-Young inequality on each cube.  We omit the details.} $f$ on $B(x_0,R)$, where $N_{1/R}(S)$ is the $1/R$-neighborhood of $S$.  To see how \eqref{eq:local-rest} implies \eqref{eq:local-rest-blur}, translate the surface $S$ by $O(1/R)$ and then average \eqref{eq:local-rest} over all such translations; to see the converse implication, introduce a bump function $\psi_{B(x_0,R)}$ concentrated near $B(x_0,2R)$ which equals 1 on $B(x_0,R)$, and exploit the reproducing formula $\hat f = \hat f * \hat \psi_{B(x_0,R)}$ to control $\hat f$ on $S$ in terms of $\hat f$ on neighborhoods such as $N_{1/R}(S)$, exploiting the fact that $\hat \psi_{B(x_0,R)}$ will decay rapidly away from $B(0,1/R)$.

Of course, \eqref{eq:local-rest-blur} is equivalent by duality to the estimate
\begin{equation}\label{eq:blur-adjoint}
\| G^\vee \|_{L^{p'}(B(x_0,R))} \leq C_{p,q,S,\alpha} R^{\alpha-1/q} \| G \|_{L^{q'}(N_{1/R}(S))}
\end{equation}
for all smooth functions $G$ supported on $N_{1/R}(S)$.  By another application of the uncertainty principle (similar to the transference principle of Marcinkiewicz and Zygmund), this estimate is also equivalent to the discrete version
\begin{equation}\label{eq:discrete-adjoint}
\| \sum_{\xi \in \Lambda} g(\xi) e^{2\pi i x \cdot \xi} \|_{L^{p'}(B(x_0,R))} \leq C_{p,q,S,\alpha} R^{\alpha+(n-1)/q} \| g \|_{l^{q'}(\Lambda)}
\end{equation}
where $\Lambda$ is any maximal $1/R$-separated subset of $S$, and $g$ is any (discrete) function on $\Lambda$.

From the formulation \eqref{eq:local-rest-blur} and Plancherel's theorem, we immediately obtain the local restriction estimate $R_S(2 \to 2; 1/2)$ for smooth compact hypersurfaces $S$; this estimate can also be obtained from the Agmon-H\"ormander estimate or from the frequency-localized version of the Sobolev trace lemma.

To convert local restriction estimates into global ones, the key tool used is the decay of the Fourier transform $(d\sigma)^\vee$.  Indeed, suppose we have a decay estimate of the form
$$ |(d\sigma)^{\vee}(x)| \leq C/(1+|x|)^\rho$$
for some $\rho > 0$.  Then the contributions to \eqref{eq:rpq} arising from widely separated portions of space will be almost orthogonal.  For instance, suppose $R \geq 1$ and $B(x_0,R)$ and $B(x_1,R)$ are two balls which are separated by at least a distance of $R$.  Then if $f_0$ and $f_1$ are supported on $B(x_0,R)$ and $B(x_1,R)$ respectively, the Fourier transforms $\hat f_0|_S$ and $\hat f_1|_S$ will be almost orthogonal on $S$:
\begin{equation}\label{eq:decay}
\begin{split}
 |\langle \hat f_0|_S, \hat f_1|_S \rangle_{L^2(S; d\sigma)}| &= |\langle \hat f_0 d\sigma, \hat f_1 \rangle_{L^2(\R^n)}|\\
 &= |\langle f_0 * (d\sigma)^\vee, f_1 \rangle| \leq C R^{-\rho} \|f_0\|_{L^1(B(x_0,R))} \| f_1\|_{L^1(B(x_1,R))},
\end{split}
\end{equation}
since the convolution kernel $(d\sigma)^\vee$ has magnitude $O(R^{-\rho})$ when applied to differences of points in $B(x_0,R)$ and points in $B(x_1,R)$.  This almost orthogonality asserts in some sense that distant balls do not interact much with each other, and so will allow us to reduce a global restriction estimate to a local one.

One heuristic way to view \eqref{eq:decay} is as follows.  This estimate is in some sense a ``bilinear'' version of the (false) estimate
\begin{equation}\label{eq:false}
 \| \hat f_0 \|_{L^2(S; d\sigma)} \leq C R^{-\rho/2} \| f_0 \|_{L^1(B(x_0,R))};
\end{equation}
this estimate is of course not true since the limit of the estimate as $R \to \infty$ is absurd, nevertheless it is ``virtually'' true in the sense that it implies the true estimate \eqref{eq:decay} by Cauchy-Schwarz.  Note that \eqref{eq:false} is just the (false) local restriction estimate $R_S(1 \to 2; -\rho/2)$.  While this estimate is not true, it is true for certain interpolation purposes; for instance, by combining it with the Agmon-H\"ormander estimate $R_S(2 \to 2; 1/2)$, one can obtain the \emph{Tomas-Stein estimate} $R_S(\frac{2(\rho+1)}{\rho+2} \to 2)$, or more generally $R_S(p \to 2)$ for all $p \leq \frac{2(\rho+1)}{\rho+2}$.  This heuristic argument can be made rigorous by using orthogonality arguments such as the $TT^*$ method; see \cite{tomas:restrict}, \cite{strichartz:restrictionquadratic}, or \cite{stein:large}.  In the particular cases of the sphere and paraboloid, the Tomas-Stein estimate yields $R_S(\frac{2(n+1)}{n+3} \to 2)$; for the cone, it yields $R_S(\frac{2n}{n-2} \to 2)$.  Note that this is consistent with the numerology supplied by the Knapp example from the previous section.

The Tomas-Stein argument uses orthogonality on $L^2(S; d\sigma)$, and at first glance it thus seems that it can only be applied to obtain restriction theorems $R_S(p \to q)$ when $q=2$.  However, it was observed by Bourgain \cite{borg:kakeya}, \cite{borg:stein} that the same type of orthogonality arguments, exploiting the decay of the Fourier transform of $d\sigma$, can also be used for restriction theorems which are not $L^2$-based, albeit with some inefficiencies due to the use of non-$L^2$ orthogonality estimates.  These ideas were then extended in \cite{vargas:restrict}, \cite{borg:cone}, \cite{tao:boch-rest}, \cite{tvv:bilinear}, \cite{tv:cone1}; we cite two sample results below.

\begin{theorem}\cite{borg:kakeya}, \cite{borg:stein}, \cite{vargas:restrict}, \cite{tvv:bilinear}, \cite{tv:cone1} Let $\rho$ be as above.  If $R^*_S(p \to q; \alpha)$ holds for some $\rho + 1 > \alpha q$, then we have $R^*(\tilde p \to \tilde q)$ whenever
$$ \tilde q > 2 + \frac{q}{\rho + 1 - \alpha q}; \quad \frac{\tilde p}{\tilde q} < 1 + \frac{q}{p(\rho + 1 - \alpha q)}.
$$
\end{theorem}

\begin{theorem}\label{thm:eps-remove}\cite{tao:boch-rest}, \cite{tao:weak2} Let $\rho$ be above.  If $R_S(p \to p; \alpha)$ holds for some $p < 2$ and $0 < \alpha \ll 1$, then we have $R_S(p \to q)$ whenever
$$ \frac{1}{q} > \frac{1}{p} + \frac{C_\rho}{\log(1/\alpha)}.$$
\end{theorem}

The second theorem in particular has the following consequence: if $R_S(p \to p; \eps)$ is true for all $\eps > 0$, then $R_S(p \to p-\eps)$ is also true for every $\eps > 0$.  (The converse statement follows easily from interpolation).  Thus we can convert a local estimate with epsilon losses to a global estimate, where the epsilon loss has now been transferred to the exponents.  This type of ``epsilon-removal lemma'' is common in this theory, see \cite{borg:cone}, \cite{tv:cone1}, \cite{tao:weak2} for some more examples.

The above results are probably not optimal, however they do emphasize the point that one can study global restriction estimates via their local counterparts.

\section{Bilinear restriction estimates}\label{sec:bilinear}

We now turn to another idea in the development of restriction theory - that of passing from the linear restriction and extension estimates to bilinear analogues.

The original motivation of this theory was the ``$L^4$'' or ``bi-orthogonality'' theory developed in such places as \cite{feff:note}, \cite{cordoba:sieve}, \cite{carl:disc}, \cite{carbery:maximal-bochner}, \cite{mock:cone}.  The basic idea is that expressions such as $\| (F d\sigma)^\vee \|_{L^{p'}(\R^n)}$ can be calculated very explicitly when $p'$ is an even integer, and especially when $p'$ is equal to 4.  Indeed, we have by Plancherel's theorem that
$$ \| (F d\sigma)^\vee \|_{L^4(\R^n)} = \| (F d\sigma)^\vee (F d\sigma)^\vee \|_{L^2(\R^n)}^{1/2} = \| F d\sigma * F d\sigma \|_{L^2(\R^n)}^{1/2}.$$
Thus one can reduce a restriction estimate such as $R_S^*(q' \to 4)$ to an estimate of the form
$$ \| F d\sigma * F d\sigma \|_{L^2(\R^n)} \leq C_q \| F \|_{L^{q'}(S; d\sigma)}^2;$$
the point here is that there is no oscillation in this estimate (since there is no Fourier transform), and this estimate can be proven or disproven by more direct methods.  For instance when $S$ is the circle in $\R^2$, there is a logarithmic divergence in the above estimate, since $d\sigma * d\sigma$ blows up like $1/|x|^{1/2}$ on the circle $\{ x \in \R^2: |x| = 2\}$ of radius 2, however by introducing the localizing parameter $R$ one can easily prove the modified estimate
\begin{equation}\label{eq:gg}
 \| G * G \|_{L^2(\R^n)} \leq C_q (\log R)^{1/2} R^{-3/2} \| G \|_{L^4(N_{1/R}(S))}^2,
\end{equation}
for all $R \geq 1$ and all $G$ supported on $N_{1/R}(S)$; comparing this with \eqref{eq:blur-adjoint} we obtain the local restriction estimate $R_S^*(4 \to 4, \eps)$ for any $\eps > 0$, which (by use of epsilon-removal lemmas such as Theorem \ref{thm:eps-remove}) proves the optimal range of restriction estimates for the circle (first proven by Zygmund \cite{zygmund}, by a more direct argument).

Similar arguments also give the optimal restriction theory for the cone in three dimensions, see \cite{barcelo}.  At first glance, this theory seems to be restricted to $L^4$, since it relies on Plancherel's theorem.  However, one can partially extend these ideas to other exponents $L^{p'}$.  The main point is that the linear estimate
$$ \| (F d\sigma)^\vee \|_{L^{p'}(\R^n)} \leq C_{p,q,S} \| F \|_{L^{q'}(S; d\sigma)}$$
is equivalent, via squaring, to the quadratic estimate
$$ \| (F d\sigma)^\vee (F d\sigma)^\vee \|_{L^{p'/2}(\R^n)} \leq C_{p,q,S} \| F \|_{L^{q'}(S; d\sigma)} \| F \|_{L^{q'}(S; d\sigma)}$$
which we can depolarize as the bilinear estimate
\begin{equation}\label{eq:bilinear}
 \| (F_1 d\sigma)^\vee (F_2 d\sigma)^\vee \|_{L^{p'/2}(\R^n)} \leq C_{p,q,S} \| F_1 \|_{L^{q'}(S; d\sigma)} \| F_2 \|_{L^{q'}(S; d\sigma)}.
\end{equation}
In such an estimate, the worst case typically occurs when $F_1$ and $F_2$ are both concentrated in the same small ``cap'' in $S$; this is what happens in the Knapp example, for instance.  

The strategy of the bilinear approach to restriction theory is to rewrite the linear estimate \eqref{eq:rpq-dual} as the bilinear estimate \eqref{eq:bilinear}, which in turn is a special case of a more general estimate of the form
\begin{equation}\label{eq:bilinear-gen}
\| (F_1 d\sigma_1)^\vee (F_2 d\sigma_2)^\vee \|_{L^{p'/2}(\R^n)} \leq C_{p,q, S_1, S_2} \| F_1 \|_{L^{q'}(S_1; d\sigma_1)} \| F_2 \|_{L^{q'}(S_2; d\sigma_2)},
\end{equation}
for arbitrary \emph{pairs} of smooth compact hypersurfaces $S_1$, $S_2$ with surface measures $d\sigma_1$, $d\sigma_2$ respectively, and all smooth $F_1$, $F_2$ supported on $S_1$ and $S_2$.  We let $R^*_{S_1,S_2}(q' \times q' \to p'/2)$ denote the statement that the estimate \eqref{eq:bilinear-gen} holds.  Then by the above discussion, $R^*_S(q' \to p')$ is equivalent to $R^*_{S,S}(q' \times q' \to p'/2)$.  Thus linear restriction estimates are special cases of bilinear estimates.  However, there are bilinear estimates that cannot be derived directly from linear ones.  For instance, let $S_1 := \{ (\xi_1,0): \xi_1 \in \R \}$ and $S_2 := \{ (0,\xi_2): \xi_2 \in \R \}$ denote the $x$ and $y$ axes in $\R^2$.  Then we have $(F_1 d\sigma_1)^\vee(x,y) = \check F_1(x)$ and $(F_2 d\sigma_2)^\vee(x,y) = \check F_2(y)$, and so there are no global restriction estimates of the form $R^*_{S_1}(q' \to p')$ or $R^*_{S_2}(q' \to p')$ unless $p' = \infty$, since the Fourier transforms do not decay at infinity.  However, since
$$(F_1 d\sigma_1)^\vee (F_2 d\sigma_2)^\vee(x,y) = \check F_1(x) \check F_2(y),$$
we see from the one-dimensional Plancherel theorem that we have the bilinear restriction estimate $R^*_{S_1, S_2}(2 \times 2 \to 2)$.  Note however that the symmetrized analogues $R^*_{S_1,S_1}(2 \times 2 \to 2)$ and $R^*_{S_2,S_2}(2 \times 2 \to 2)$ are false.  Thus the bilinear estimate exploits the \emph{transversality} of $S_1$ and $S_2$.

A higher-dimensional analogue of this estimate is known: if $S_1$ and $S_2$ are two smooth compact hypersurfaces which are \emph{transverse} in the sense that the set of unit normals of $S_1$ are separated by some non-zero distance from the set of unit normals of $S_2$, then we have $R^*_{S_1, S_2}(2 \times 2 \to 2)$.  This can be easily seen by using Plancherel to convert the bilinear restriction estimate to a bilinear convolution estimate
$$
\| (F_1 d\sigma_1) * (F_2 d\sigma_2) \|_{L^{2}(\R^n)} \leq C_{p,q, S_1, S_2} \| F_1 \|_{L^{2}(S_1; d\sigma_1)} \| F_2 \|_{L^{2}(S_2; d\sigma_2)},$$
and then using the Cauchy-Schwarz estimate
$$
\| (F_1 d\sigma_1) * (F_2 d\sigma_2) \|_{L^{2}(\R^n)}
\leq \| (|F_1|^2 d\sigma_1) * (|F_2|^2 d\sigma_2) \|_{L^1(\R^n)}
\| d\sigma_1) * d\sigma_2 \|_{L^\infty(\R^n)}$$
and using transversality to bound the second factor.  Generalizations of these ``bilinear $L^2$'' estimates have arisen in recent work in non-linear evolution equations (starting with the work of Bourgain \cite{borg:xsb} and Klainerman-Machedon \cite{klainerman:nulllocal} and continued by numerous authors, see for instance \cite{kpv:kdv}) and are especially useful for handling non-linearities which contain derivatives arranged to create a ``null form'', but we will not pursue this matter here, and refer the reader instead to \cite{ginibre:survey}, \cite{damiano:null}, \cite{tao:xsb}.  There has been also some work in generalizing these bilinear estimates to weighted settings, see \cite{bcc}.

Now let $S_1$ and $S_2$ be two compact transverse subsets of the sphere or paraboloid.  Of course, any restriction theorem $R^*_S(q' \to p')$ for the sphere or paraboloid will imply a bilinear restriction theorem $R^*_{S_1,S_2}(q' \times q' \to p'/2)$ for the pair $S_1$, $S_2$.  However, the transversality allows us to prove more estimates in the bilinear setting than the linear one; we have already seen the estimate $R^*_{S_1,S_2}(2 \times 2 \to 2)$, whereas the linear restriction estimate $R^*_S(2 \to 4)$ is only true in three and higher dimensions.  One reason for this is that there is no exact analogue of the Knapp example in the transverse bilinear setting.  Indeed, the best necessary conditions known on $R^*_{S_1, S_2}(q' \times q' \to p'/2)$ are that
\begin{equation}\label{eq:bil-nec}
 p > \frac{2n}{n+1}; \quad \frac{n+2}{p'} + \frac{n}{q'} \leq n; \quad \frac{n+2}{p'} + \frac{n-2}{q'} \leq n-1;
\end{equation}
see \cite{tvv:bilinear} or \cite{damiano:null}, where one develops bilinear analogues of the Knapp examples.  This is somewhat less stringent than the corresponding conditions
\begin{equation}\label{eq:knapp}
 p > \frac{2n}{n+1}; \quad \frac{n+1}{p'} + \frac{n-1}{q'} \leq n-1
\end{equation}
for the linear problem $R^*_S(q' \to p')$.  The \emph{bilinear restriction conjecture} asserts that the necessary conditions \eqref{eq:bil-nec} are in fact sharp.  This conjecture is still open except in dimension two, but recently it has been shown that (up to endpoints) it is equivalent to the usual restriction conjecture for the sphere and paraboloid.

Up until now, we have viewed bilinear restriction estimates as being more complex generalizations of linear restriction estimates, which seems to offer no incentive to study the bilinear estimates until the linear ones are settled.  However, it turns out that one can use the bilinear estimates to go back and deduce new linear estimates, and indeed all the recent progress on the restriction problem has been obtained in this manner.  The key observation is that  one can perform a Whitney decomposition of the product manifold $S \times S$ around the diagonal $\Delta := \{ (\xi,\xi): \xi \in S\}$ so that $S \times S \backslash \Delta$ decomposes as the disjoint union of sets of the form $S_1 \times S_2$, where $S_1, S_2$ are disjoint subsets of $S$ whose separation is comparable to their diameter.  This allows one to obtain bilinear restriction estimates of the form $R^*_{S,S}(p' \times p' \to q'/2)$ (and hence $R^*_S(p' \to q')$) from estimates of the form $R^*_{S_1,S_2}(p' \times p' \to q'/2)$, using some rescaling and orthogonality estimates to sum up (and discarding the diagonal $\Delta$, which is of measure zero); see \cite{tvv:bilinear} for more details.  Of course one cannot hope to have an unconditional implication of the form $R^*_{S_1,S_2}(p' \times p' \to q'/2) \implies R^*_S(p' \to q')$, since the necessary conditions \eqref{eq:bil-nec} for the former are weaker than those \eqref{eq:knapp} for the latter; however, we can do the next best thing:

\begin{theorem}\label{thm:tvv}\cite{tvv:bilinear}  Let $p$, $q$ obey the necessary conditions \eqref{eq:knapp}, and suppose that $R^*_{S_1,S_2}(\tilde p' \times \tilde p' \to \tilde q')$ is true for all $(\tilde p, \tilde q)$ in an open neighborhood of $(p,q)$, and some pair $S_1, S_2$ of compact transverse subsets of the paraboloid.  Then $R^*_S(p' \to q')$ is true.
\end{theorem}

A similar result is true for the sphere, except that one must make $S_1$ and $S_2$ subsets of a certain parabolically rescaled version of the sphere; see \cite{tvv:bilinear} for more details.

The above theorem (and ones like it) allow one to pass freely back and forth between linear and (transverse) bilinear restriction estimates.  For instance, this theorem can be used to provide an alternative proof of Zygmund's estimate (which asserts in particular that $R^*_S(4 \to 4+\eps)$ when $S$ is the unit circle) from the basic estimate $R^*_{S_1,S_2}(2 \times 2 \to 2)$ for transverse sets.  Although the bilinear estimates appear more complicated, they are in fact easier to analyze because they consist purely of \emph{transverse interactions}, excluding the \emph{parallel interactions} which often cause the most trouble (cf. the Knapp example).

One can of course formulate local bilinear restriction estimates $R^*_{S_1,S_2}(p' \times p' \to q'; \alpha)$, which assert that
$$
\| (F_1 d\sigma_1)^\vee (F_2 d\sigma_2)^\vee \|_{L^{p'/2}(B(x_0,R))} \leq C_{p,q, S_1, S_2, \alpha} R^\alpha \| F_1 \|_{L^{q'}(S_1; d\sigma_1)} \| F_2 \|_{L^{q'}(S_2; d\sigma_2)}.
$$
One can of course reformulate these estimates using the uncertainty principle in a similar way to before, though some reformulations are not available because the notion of dualizing a bilinear estimate becomes difficult to use.  There are also bilinear ``epsilon-removal'' lemmas available; for instance, we have

\begin{theorem}\label{thm:eps}\cite{borg:cone}, \cite{tv:cone1}  Let $S_1$, $S_2$ be compact surfaces obeying some decay estimate 
$$ |(d\sigma_1)^\vee(x)|, |(d\sigma_2)^\vee(x)| \leq C/(1+|x|)^\rho$$
for some $\rho > 0$.  Suppose we have the local bilinear restriction estimate $R^*_{S_1,S_2}(2 \times 2 \to q, \eps)$ for all $\eps > 0$.  Then we have the global bilinear restriction estimate $R^*_{S_1,S_2}(2 \times 2 \to q+\eps)$ for all $\eps > 0$.
\end{theorem}

More quantitative versions of this estimate have been proven, see e.g. \cite{tv:cone1}, Lemma 2.4.  See also \cite{krt} for a more PDE-based approach to this epsilon-removal lemma.

The bilinear estimate $R^*_{S_1,S_2}(2 \times 2 \to 2)$ holds for all surfaces $S_1$, $S_2$ which are transverse.  If both $S_1$ and $S_2$ are flat, then this estimate is sharp; however one can improve this estimate slightly when $S_1$ and $S_2$ have some curvature.  For instance, if $S_1$ and $S_2$ are transverse subsets of the paraboloid in $\R^n$, then we have $R^*_{S_1,S_2}(p \times p \to 2)$ for all $p \geq \frac{4n}{3n-2}$, see \cite{tvv:bilinear}.  To see why we should gain over the $p=2$ estimate, consider the following.  Using Plancherel, we can rewrite $R^*_{S_1,S_2}(2 \times 2 \to 2)$ as
the bilinear convolution estimate
$$ \| (F_1 d\sigma_1) * (F_2 d\sigma_2) \|_{L^2(\R^n)}^2 \leq C \| F_1 \|_{L^2(S_1)}^2 \| F_2 \|_{L^2(S_2)}^2.$$
Let us suppose for the moment that we are in a model case, where the $F_j$, $j=1,2$ are characteristic functions, for the sets $\{ (\underline{\xi_j}, \frac{1}{2} |\underline{\xi_j}|^2): \underline{\xi_j} \in \Omega_j \}$ for some disjoint bounded open subsets $\Omega_1$, $\Omega_2$ of $\R^n$; the right-hand side is thus $C |\Omega_1| |\Omega_2|$.  Then by discarding some Jacobian factors (which are harmless due to the transversality), the left-hand side is essentially the volume of the $3n-1$-dimensional set
$$ \{ (\underline{\xi}_1, \underline{\xi}_2, \underline{\xi}_3, \underline{\xi}_4) \in \Omega_1 \times \Omega_2 \times \Omega_1 \times \Omega_2:
 \underline{\xi}_1+\underline{\xi}_2 = \underline{\xi}_3 + \underline{\xi}_4; \quad |\underline{\xi}_1|^2 + |\underline{\xi}_2|^2 = |\underline{\xi}_3|^2 + |\underline{\xi}_4|^2\}.$$
The two constraints imply that $\underline{\xi}_1$, $\underline{\xi}_2$ and $\underline{\xi}_3$, $\underline{\xi}_4$ from the opposing diagonals of a rectangle.  In particular, $\underline{\xi}_3$ lies on the hyperplane $\pi(\underline{\xi}_1, \underline{\xi}_4)$ containing $\underline{\xi}_1$ and orthogonal to $\underline{\xi}_4 - \underline{\xi}_1$, and then $\underline{\xi}_2$ can be recovered from the other three frequencies by the formula $\underline{\xi}_2 = \underline{\xi}_3 + \underline{\xi}_4 - \underline{\xi}_1$.  Thus, by Fubini's theorem, the volume of the above set is bounded above by
$$ \int_{\Omega_1} \int_{\Omega_2} | \Omega_2 \cap \pi(\underline{\xi}_1, \underline{\xi}_4) |\ d\underline{\xi}_1 d\underline{\xi}_4$$
(discarding the constraint $\xi_2 \in \Omega_2$).  Since $\Omega_2$ is bounded, we may make the very crude estimate
\begin{equation}\label{eq:crude}
 | \Omega_2 \cap \pi(\underline{\xi}_1, \underline{\xi}_4) | \leq C,
\end{equation}
from which the desired bound of $C |\Omega_1| |\Omega_2|$ follows.

The estimate \eqref{eq:crude} can of course be improved when $\Omega_2$ is small, using the standard $L^p$ bounds for the Radon transform.  This is made rigorous in \cite{vargas:restrict}, \cite{vargas:2}, \cite{tvv:bilinear}, culminating in the above-mentioned bilinear restriction estimate $R^*_{S_1,S_2}(p \times p \to 2)$ for all $p \geq \frac{4n}{3n-2}$.  This issue of exploiting the possible gain over \eqref{eq:crude} also arises in some recent developments \cite{tao:parabola} in bilinear restriction theory, which we shall return to later.

The latter two conditions of \eqref{eq:knapp} meet when $q'=2$, when they assert that $R^*_{S_1, S_2}(2 \times 2 \to q)$ for all $q \geq \frac{2(n+2)}{n}$.  This was first conjectured by Machedon and Klainerman for both the paraboloid and cone.  Despite the original restriction conjecture remaining open, this conjecture has been completely solved for the cone and solved except for an endpoint for the paraboloid; see Figures 3, 4.  We shall discuss this recent progress in the next few sections.

\begin{figure}\label{fig:cone-bil}
\begin{tabular}{|l|l|l|} \hline
Dimension & Range of $q$ & \\
\hline
$n \geq 2$  & $q \geq 2$ & Plancherel + Cauchy-Schwarz\\ 
$n \geq 3$  & $q \geq \frac{n}{n-2}$ & Strichartz 1977 
\cite{strichartz:restrictionquadratic} \\
$n = 2$  & $q \geq 2 - \frac{13}{2408}$ & Bourgain 1995 \cite{borg:cone}\\ 

$n = 2$  & $q \geq 2 - \frac{8}{121}$ & Tao, Vargas 2000 \cite{tv:cone1}\\ 
$n \geq 2$ & $q > \frac{2(n+2)}{n}$ & Wolff, 2000 \cite{wolff:cone}\\
$n \geq 2$ & $q \geq \frac{2(n+2)}{n}$ & Tao, 2001 \cite{tao:cone} (best possible)\\
\hline
\end{tabular}
\caption{Known results on the bilinear restriction problem $R_{S_1 \times S_2}(2 \times 2 \to q)$, for transverse compact subsets of the cone. }
\end{figure}

\begin{figure}\label{fig:parab-bil}
\begin{tabular}{|l|l|l|} \hline
Dimension & Range of $q$ & \\
\hline
$n \geq 2$  & $q \geq 2$ & Plancherel + Cauchy-Schwarz\\ 
$n \geq 3$  & $q \geq \frac{n+1}{n-1}$ & Strichartz 1977 
\cite{strichartz:restrictionquadratic} \\
$n = 2$  & $q \geq 2 - \frac{5}{69}$ & Tao, Vargas, Vega 1998 \cite{tvv:bilinear}\\ 
$n = 2$  & $q \geq 2 - \frac{2}{17}$ & Tao, Vargas 2000 \cite{tv:cone1}\\ 

$n \geq 2$ & $q > \frac{2(n+2)}{n}$ & Tao, 2003 \cite{tao:parabola}\\
$n \geq 2$ & $q \geq \frac{2(n+2)}{n}$ & (conjectured)\\
\hline
\end{tabular}
\caption{Known results on the bilinear restriction problem $R_{S_1 \times S_2}(2 \times 2 \to q)$, for transverse compact subsets of the paraboloid. }
\end{figure}

\section{The wave packet decomposition}\label{sec:wave-packet}

We have discussed two of the tools in the modern theory of restriction estimates: the reduction to local estimates, and the reduction to bilinear estimates.  We now turn to a third key technique: the introduction of wave packets.

For sake of illustration, suppose we wish to prove the local restriction estimate $R_S(p \to 1; \alpha)$ where $S$ is the sphere \eqref{eq:sphere}; the exponent $1$ can of course be changed, but this does not significantly alter the argument sketched below (except that the estimates on certain coefficients $c_T$ will change).  We use the formulation \eqref{eq:blur-adjoint}, fixing $x_0 = 0$, thus we have to prove
$$\| G^\vee \|_{L^{p'}(B(0,R))} \lesssim R^{\alpha}$$
for $R \geq 1$ all bounded functions $G$ on $N_{1/R}(S)$.  Henceforth we use $X \lesssim Y$ or $X = O(Y)$ to denote the estimate $X \leq CY$ for some $C$ depending on such parameters as $n$, $p$, $\alpha$, while we call a function ``bounded'' when it has an $L^\infty$ norm of $O(1)$.

Fix $R$ and $G$, and observe that the annular region $N_{1/R}(S)$ can be divided into $\sim R^{(n-1)/2}$ finitely overlapping disks $\kappa$ of width $\sim 1/\sqrt{R}$ and thickness $1/R$.  If $G$ is a function on $N_{1/R}(S)$, we can thus use a partition of unity to divide $G = \sum_\kappa G_\kappa$, where each $G_\kappa$ is a bounded function supported on one of these disks $G_\kappa$.  Our task is thus to show that
$$\| \sum_\kappa G_\kappa^\vee \|_{L^{p'}(B(0,R))} \lesssim R^{\alpha}.$$
The question then arises as to what $G_\kappa^\vee$ looks like.  We first consider some examples.  Suppose that the disk $G_\kappa$ is centered at a point $\omega_\kappa \in S^{n-1}$, which by the geometry of the sphere implies that $\omega_\kappa$ is also essentially the normal to the disk $\kappa$.  If $G_\kappa$ is a bump function adapted to $\kappa$, then by duality $G_\kappa^\vee$ would be concentrated on the $R \times R^{1/2}$ tube
$$ T_{\kappa,0} := \{ x \in B(0,R): \pi_{\omega_\kappa^\perp} x = O(R^{1/2}) \},$$
where $\pi_{\omega_\kappa^\perp}$ is the orthogonal projection onto the hyperplane $\omega_\kappa^\perp := \{ x \in \R^n: x \cdot \omega_\kappa = 0\}$.  Indeed, since $\kappa$ has volume roughly $R^{-(n+1)/2}$, we would expect $G_\kappa^\vee$ to equal a function $\psi_{T_{\kappa,0}}$ of the form
\begin{equation}\label{eq:wavepacket}
\psi_{T_{\kappa,0}}(x) = R^{-(n+1)/2} e^{2\pi i \omega_\kappa \cdot x} \phi_{T_{\kappa,0}},
\end{equation}
where $\phi_{T_{\kappa,0}}$ is a Schwartz function adapted to the tube $T_{\kappa,0}$ which has size $O(1)$ on this tube and is rapidly decreasing away from this tube.  We call the function $\psi_{T_{\kappa,0}}$ a \emph{wave packet} adapted to the tube $T_{\kappa,0}$; this object has already essentially come up in the discussion of the Knapp example.

What happens when $G_\kappa$ is not a bump function adapted to $\kappa$?  First suppose that $G_\kappa$ is a modulated bump function, more precisely suppose
$$ G_\kappa(\xi) = e^{-2\pi i x_0 \cdot \xi} \tilde G_\kappa(\xi)$$
where $\tilde G_\kappa$ is a bump function adapted to $\kappa$, and $x_0$ is an element of the hyperplane $\omega_\kappa^\perp$.  Then by the above discussion, $G_\kappa^\vee$ will be concentrated on the $R \times R^{1/2}$ tube
$$ T_{\kappa,x_0} := T_{\kappa,0} + x_0,$$
indeed we have
$$
G_\kappa^\vee(x) = \psi_{T_{\kappa,x_0}} := R^{-(n+1)/2} e^{2\pi i \omega_\kappa \cdot x} \phi_{T_{\kappa,x_0}}$$
for some Schwartz function $\psi_{T_{\kappa,x_0}}$ adapted to $T_{\kappa,x_0}$.
(One could also modulate $G_\kappa$ in the direction parallel to $\omega_\kappa$ instead of in the perpendicular directions, but this either has a negligible effect on the Fourier transform on the ball $B(0,R)$, or else makes the Fourier transform much smaller, depending on how much modulation is applied).

Thus one can make $G_\kappa^\vee$ resemble a wave packet $\psi_T$ for any tube $T$ oriented in the direction $\omega_\kappa$.  In the general situation, where $G_\kappa$ is a bounded function on $\kappa$, then one can perform a Fourier series decomposition in the directions perpendicular to $\omega_\kappa$ to essentially decompose $G_\kappa$ as an $l^2$-average of modulated bump functions.  (The behavior in the direction parallel to $\omega_\kappa$, which only extends for a distance $O(1/R)$ is essentially irrelevant, thanks to the localization of physical space to $B(0,R)$ and the uncertainty principle).  Thus we can write
$$ G_\kappa^\vee = \sum_{T // \omega_\kappa} c_T \psi_T,$$
where $T$ ranges over a finitely overlapping collection of $R \times \sqrt{R}$ tubes in $B(0,R)$ oriented in the direction $\omega_\kappa$, $\psi_T$ is a wave packet adapted to $T$, and $c_T$ is a collection of scalars with the $L^2$ normalization condition $\sum_{T // \omega_\kappa} |c_T|^2 \lesssim 1$.
One can then expand the original Fourier transform $G^\vee$ as
$$ G^\vee = \sum_T c_T \psi_T$$
where $T$ now ranges over a separated\footnote{This means that any two tubes $T$, $T'$ in this collection either have directions differing by at least $1/R^{1/2}$, or are parallel and are separated spatially by at least $R^{1/2}$.} collection of tubes in $B(0,R)$, and $\omega_T$ denotes the direction of $T$.

This heuristic decomposition is an example of what is known as the \emph{wave packet decomposition} of $G^\vee$.  Versions of this decomposition in the context of the restriction problem (or the closely related Bochner-Riesz problem) first appeared in \cite{cordoba:covering}, \cite{cordoba:sieve}, \cite{feff:note}, \cite{feff:ball}, and was then later developed in \cite{borg:kakeya}, \cite{borg:stein}, \cite{vargas:restrict}, \cite{vargas:2}, \cite{tvv:bilinear}, \cite{tv:cone1}, \cite{wolff:cone}, \cite{tao:cone}; this method also can be applied to related problems such as local smoothing or Bochner-Riesz, see for instance \cite{wolff:distance}, \cite{wolff:smsub}.  The wave packet decomposition reduces the study of restriction estimates to that of proving estimates on the linear superpositions of wave packets
\begin{equation}\label{eq:osc}
 \| \sum_T c_T \psi_T \|_{L^{p'}(B(0,R))}.
\end{equation}
Note that the wave packets $\psi_T$ have two main features; one at coarse scales $\gg \sqrt{R}$ and one at fine scales $\ll \sqrt{R}$.  At coarse scales, the wave packet is localized to a relatively thin tube of width $\sqrt{R}$.  At fine scales, the wave packet oscillates at a fixed frequency $\omega_T$.  Note that the coarse scale behavior and fine scale behavior are linked, because the direction of the tube at coarse scales is exactly the same as the frequency of the oscillation at fine scales.  The issue is then how to co-ordinate these two aspects - localization at coarse scales, and oscillations at fine scales - of wave packets in order to estimate \eqref{eq:osc} efficiently.

The first strategy for estimating these superpositions of wave packets is due to C\'ordoba \cite{cordoba:covering}, \cite{cordoba:sieve}, in which the idea is to estimate the oscillatory sum by the associated square function
\begin{equation}\label{eq:sq}
 \| (\sum_T |c_T \psi_T|^2)^{1/2} \|_{L^{p'}(B(0,R))}.
\end{equation}
The point of doing so is that all the fine-scale oscillation has been removed from this problem, leaving only the coarse scale localizations to tubes.  There is still of course the problem of estimating this non-oscillatory square function; this problem is essentially equivalent\footnote{Conversely, one must resolve the Kakeya conjectures in order to fully resolve the restriction problem, because one can use randomization arguments to show that any bound on \eqref{eq:osc} implies a comparable bound on \eqref{eq:sq}.  See e.g. \cite{bcss}.} to the problem of estimating the \emph{Kakeya maximal function}, which is another important problem in harmonic analysis, but one which we will not discuss in detail here.  (See however \cite{wolff:kakeya}, \cite{Bo}, \cite{tao:elesc}).

Now we discuss how to estimate the oscillatory sum \eqref{eq:osc} by the square function \eqref{eq:sq}.  When $p' = 2$, or when $n=2$ and $p' = 4$, one can bound the former by the latter by direct orthogonality (or bi-orthogonality) arguments, however these arguments do not work for other values of $p'$.  Nevertheless, it was observed by Bourgain \cite{borg:kakeya}, \cite{borg:stein} that one can still obtain some control of \eqref{eq:osc} by \eqref{eq:sq} in these cases, but with a loss of some powers of $R$.  The idea is to break the ball $B(0,R)$ up into cubes $q$ of size $\sqrt{R}$.  On such ``fine-scale'' cubes, a wave packet $\psi_T = R^{-(n-1)/2} e^{2\pi i \omega_T \cdot x} \phi_T$ has essentially constant magnitude; to (over-)simplify the discussion, let us suppose that $\phi_T$ is equal to 1 on $q$ if $q \subset T$ and $\phi_T$ vanishes on $q$ otherwise.  Then the portion of \eqref{eq:osc} coming from $q$ is
$$ R^{-(n-1)/2} \| (\sum_{T: T \supset q} c_T e^{2\pi i \omega_T \cdot x} \|_{L^{p'}(B(0,R))}$$
while the corresponding portion of \eqref{eq:sq} is essentially
$$ R^{-(n-1)/2} R^{n/2p'} (\sum_{T: T \supset q} |c_T|^2)^{1/2}.$$
One can then control the former expression by the latter using discrete restriction estimates\footnote{It is intriguing that one uses local restriction estimates at scale $\sqrt{R}$, together with some Kakeya information, to obtain local restriction estimates at scale $R$.  This suggests a possible ``bootstrap'' approach where one could continually improve restriction estimates via iteration.  Some partial iteration methods to this effect can be found in \cite{borg:cone}, \cite{tvv:bilinear}, \cite{tv:cone1}; another example of this idea occurs in the induction-on-scales approach discussed in the next section.} of the type \eqref{eq:discrete-adjoint}, although the various powers of $R$ which accumulate when doing so do not necessarily all cancel, and so this method of estimation can cause some losses\footnote{It is conjectured that in any dimension $n \geq 2$, that one can estimate \eqref{eq:osc} by \eqref{eq:sq} in the endpoint case $p = 2n/(n+1)$, with at most an epsilon loss $R^\eps$; this, together with the so-called Kakeya maximal function conjecture, would imply the restriction conjecture.  However, it is nowhere near solved at present, except when $n=2$, and is likely to be a harder problem than the restriction problem itself.}. 

By combining these observations with some non-trivial progress on the Kakeya maximal function conjecture, Bourgain \cite{borg:kakeya}, \cite{borg:stein} was able to obtain certain improvements to the Tomas-Stein estimate (see Figure \eqref{fig:parab}).  Further progress was made by Wolff \cite{wolff:kakeya}, who improved the Kakeya estimate used in Bourgain's argument.  By introducing bilinear (or $L^4$) methods to these arguments, further improvements were obtained in \cite{vargas:restrict}, \cite{tvv:bilinear}, \cite{borg:cone} \cite{tv:cone1}; one feature of these bilinear methods is that they could now be applied to the cone as well as the sphere or paraboloid.  

These methods, however, did not obtain sharp ranges of exponents, for a variety of technical reasons.  The next breakthrough was achieved by Wolff \cite{wolff:cone}, who solved (up to endpoints) the Machedon-Klainerman conjecture for cones, by employing one additional technique - that of induction on scales, which we discuss next.

\section{Induction on scales}\label{sec:induct}

The strategy to prove a local restriction estimate at a scale $R$ in the previous section can be summed up as follows: starting with a function $G^\vee$, decompose it into wave packets supported on $\sqrt{R} \times R$ tubes.  Designating scales greater than $\sqrt{R}$ as coarse, and scales less than $\sqrt{R}$ as fine, we use oscillatory estimates such as local restriction
estimates on fine scales, and Kakeya type estimates at coarse scales, in order to obtain the desired control on $G^\vee$.

This type of argument works particularly well when $G$ is a Knapp example supported on a disk of radius $R^{-1/2}$, so that $G^\vee$ is essentially a single wave packet.  However, it becomes inefficient when $G$ is a Knapp example spread out over a wider region, e.g. a cap-type region of radius $r^{-1/2}$ for some $1 \leq r \leq R$.  Then $G^\vee$ is concentrated on a much smaller set than a single wave packet - indeed, it is (somewhat) localized to a $r \times \sqrt{r}$ tube instead of an $R \times \sqrt{R}$ tube - but the wave packet decomposition requires that one decompose $G^\vee$ as the sum of much larger objects.  This is a rather inefficient decomposition, and one which leads to significant losses in the estimates.

The difficulty here is that the wave packet decomposition is chosen in advance, instead of being adapted to the particular function $G$ being investigated.  In particular, if it turns out that $G^\vee$ is concentrating in a much smaller region, say a ball $B(x_0,r)$, then one should replace the rather coarse $R \times \sqrt{R}$ wave packet decomposition by a finer one, in this case a $r \times \sqrt{r}$ decomposition.

Of course, the difficulty is that it would be incredibly complicated to actually try to construct such an adaptive wave packet decomposition, recursively passing from coarser scales to finer scales.  Fortunately, a way out of this complexity was discovered by Wolff \cite{wolff:cone} - which is to hide all this recursive complexity in an induction hypothesis, which we now refer to as an \emph{induction on scales} argument.  Using this new idea, Wolff was able to obtain a nearly-sharp bilinear restriction estimate for the cone, namely that $R^*_{S_1,S_2}(2 \times 2 \to q)$ is true for all transverse compact subsets $S_1$, $S_2$ of the cone, and all $q > \frac{n+2}{n}$.

We now describe, rather informally, the idea of the argument; for a more rigorous presentation see \cite{wolff:cone}, \cite{tao:cone}, \cite{tao:non-endpoint}, \cite{tao:parabola}, \cite{krt}.  Suppose inductively that we already have some local estimate of the form $R^*_S(q' \to p'; \alpha)$; we will now try to use this estimate to prove a better estimate of the form $R^*_S(q' \to p'; \alpha - \eps)$ for some $\eps > 0$ depending on $\alpha$ (in what follows, the value of $\eps$ will vary from line to line).  Iterating this, we will eventually be able to obtain the estimate $R^*_S(q' \to p', \eps)$ for any $\eps > 0$, at which point we can use epsilon removal lemmas to obtain a global restriction estimate.

We still have to obtain the estimate $R^*(q' \to p'; \alpha - \eps)$ from the inductive hypothesis $R^*(q' \to p'; \alpha)$.  We first describe a somewhat oversimplified version of the main idea as follows.  We have to prove an estimate of the form
$$
\| (F d\sigma)^\vee \|_{L^{p'}(B(0,R))} \leq C_{p,q,S,\alpha} R^{\alpha-\eps} \| F \|_{L^{q'}(S)}
$$
for some $F$ on $S$ and $R \geq 1$, which we now fix.  Introduce the scale $r := R^{1-\eps}$, which is slightly smaller than $R$.  Then by the induction hypothesis $R^*(q' \to p'; \alpha)$ applied to scale $r$, we have
$$
\| (F d\sigma)^\vee \|_{L^{p'}(B(x_0,r))} \leq C_{p,q,S,\alpha} R^{\alpha-\eps} \| F \|_{L^{q'}(S)}$$
for any ball $B(x_0,r)$.  Thus we can already prove the desired estimate on smaller balls $B(x_0,r)$.  More generally, we can prove
$$
\| (F d\sigma)^\vee \|_{L^{p'}(\bigcup_j B(x_j,r))} \leq C_{p,q,S,\alpha} R^{\alpha-\eps} \| F \|_{L^{q'}(S)}$$
on any union $\bigcup_j B(x_j,r)$ of smaller balls, as long as the number of balls involved is not too large (e.g. at most $O((\log R)^C)$ for some absolute constant $C$).

As a rough first approximation, the idea of Wolff is to identify the ``bad'' balls $B(x_j,r)$ on which the function $(F d\sigma)^\vee$ ``concentrates''; the choice of these balls will of course depend on $F$.  These balls can be dealt with using the induction hypothesis, and it then remains to verify the restriction estimate on the exterior of these bad balls:
$$
\| (F d\sigma)^\vee \|_{L^{p'}(B(0,R) - \bigcup_j B(x_j,r))} \leq C_{p,q,S,\alpha} R^{\alpha-\eps} \| F \|_{L^{q'}(S)}.
$$

The above description of Wolff's argument was something of an oversimplification for two reasons; firstly, Wolff is working in the bilinear setting rather than the linear setting, and secondly the balls $B(x_j,r)$ turn out to depend not only on the original function $F$, but of the wave packet decomposition associated to $F$.  Let us ignore the first reason for the moment, and clarify the second.  On the ball $B(0,R)$, one can obtain a wave packet decomposition of the form
$$ (F d\sigma)^\vee(x) = \sum_T c_T \psi_T.$$
Because the argument of Wolff dealt with the cone, the wave packet decomposition here is slightly different from that discussed in the previous section, in two respects: firstly, the tubes $T$ are oriented on ``light rays'' normal to the cone $S$ instead of pointing in general directions, and secondly the internal structure of the wave packet $\psi_T$ is more interesting than just the product of a plane wave and a bump function, being decomposable into ``plates''.  We however will gloss over this technical issue.

For simplicity, let us suppose that the constants $c_T$ behave like a characteristic function; more precisely, there is some collection $\T$ of tubes such that $c_T = c$ for $T \in \T$ and $c_T = 0$ otherwise.  (The general case can be reduced to this case via a dyadic pigeonholing argument, which costs a relatively small factor of $\log R$).  Then we have
$$ (F d\sigma)^\vee(x) = c \sum_{T \in \T} \psi_T(x).$$
The idea now is to allow each wave packet $\psi_T$  to be able to ``exclude'' a single ball $B_T$ of the slightly smaller radius $r$.  In other words, one divides $(F d\sigma)^\vee$ into two pieces, a ``localized'' piece
$$ c \sum_{T \in \T} \psi_T(x) \chi_{B_T}(x)$$
and the ``global'' piece
$$ c \sum_{T \in \T} \psi_T(x) (1 - \chi_{B_T}(x)).$$
One then tries to control the localized piece using the induction hypothesis, and then handle the non-localized piece using the strategy of the previous section.

In the linear setting, this strategy does not quite work, because the localized pieces cannot be adequately controlled by the induction hypothesis.  However, in the bilinear setting, when one is trying to prove an estimate of the form
$$
\| (F_1 d\sigma_1)^\vee (F_2 d\sigma_2)^\vee \|_{L^{p'/2}(B(0,R))} \leq C_{p,q,S,\alpha} R^{\alpha-\eps} \| F_1 \|_{L^{q'}(S)} \| F_2 \|_{L^{q'}(S)}
$$
then one can decompose
$$ (F_j d\sigma_j)^\vee(x) = c_j \sum_{T_j \in \T_j} \psi_{T_j}(x)$$
for $j=1,2$, and allow each tube $T_j$ to exclude a single\footnote{Actually, in Wolff's argument there are $O((\log R)^C)$ such balls excluded, but this is a minor technical detail.} ball $B_{T_j}$ of radius $r$.  We can then split the bilinear expression
$$ \sum_{T_1 \in \T_1} \sum_{T_2 \in \T_2} \psi_{T_1} \psi_{T_2}$$
into a local piece
$$ \sum_{T_1 \in \T_1} \sum_{T_2 \in \T_2} \psi_{T_1} \psi_{T_2} \chi_{B_{T_1} \cap B_{T_2}}$$
(where \emph{both} tubes $T_1$ and $T_2$ are excluding $x$), and a global piece
$$ \sum_{T_1 \in \T_1} \sum_{T_2 \in \T_2}  \psi_{T_1} \psi_{T_2} (1 - \chi_{B_{T_1} \cap B_{T_2}}).$$
The local piece turns out to be easily controllable by the inductive hypothesis, so it remains to control the global piece.

The key point is to prevent too many of the tubes $T_1$ and $T_2$ from interacting with each other.  This is done by selecting the balls $B_{T_1}$, $B_{T_2}$ strategically.  Roughly speaking, for each tube $T_1$, we choose $B_{T_1}$ to be the ball which contains as many intersections of the form $T_1 \cap T_2$ as possible; the ball $B_{T_2}$ is chosen similarly.  The effect of this choice is that any point $x$ which lies in a large number of tubes in $T_1$ and in $T_2$ simultaneously, is likely to be placed primarily in the local part of the bilinear expression, and not in the global part.

With this choice of the excluding balls $B_{T_1}$, $B_{T_2}$, Wolff was able to obtain satisfactory control on the number of times tubes $T_1$ from $\T_1$ would intersect tubes $T_2$ from $\T_2$.  The key geometric observation is as follows.  Suppose that many tubes $T_1$ in $\T_1$ were going through a common point $x_0$; since the tubes $T_1$ are constrained to be oriented along light rays, these tubes must then align on a ``light cone''.

Now consider a tube $T_2$ from $\T_2$; this tube is of course transverse to all the tubes $T_1$ considered above, and furthermore is transverse to the light cone that the tubes $T_1$ lie on.  It can either pass near $x_0$, or stay far away from $x_0$.  In the first case it turns out that the joint contribution of the tubes $T_1$ and $T_2$ will largely lie in the local part of the bilinear expression and thus be manageable.  In the second case we see from transversality that the tube $T_2$ can only intersect a small number of tubes $T_1$.

Thus if there is too much intersection among tubes $T_1$ in $\T_1$, then there will be fairly sparse interection between those tubes $T_1$ and tubes $T_2$ in $\T_2$.  This geometric fact was exploited via combinatorial arguments in \cite{wolff:cone}, and when combined with some local $L^2$ arguments from \cite{mock:cone} to handle the fine scale oscillations, and the induction on scales argument, was able to obtain the near-optimal bilinear restriction theorem $R^*_{S_1, S_2}(2 \times 2 \to q)$ for $q > \frac{n+2}{n}$.  (The endpoint $q = \frac{n+2}{n}$ to the Machedon-Klainerman conjecture was then obtained in \cite{tao:cone} by refining the above argument).

\section{Adapting Wolff's argument to the paraboloid}\label{sec:parab}

The above argument of Wolff \cite{wolff:cone}, which yielded the optimal bilinear $L^2$ restriction theorem for the cone, relied on a key fact about the cone: all the tubes passing through a common point $x_0$, were restricted to lie on a hypersurface (specifically, the cone with vertex at $x_0$).  This property does not hold for the paraboloid, since in this setting the tubes can point in arbitrary directions.  Nevertheless, it is possible to recover this hypersurface property by exploiting a little more structure at fine scales, and more precisely by squeezing one ``dimension'' of gain out of \eqref{eq:crude}, thus obtaining the optimal bilinear $L^2$ restriction theorem for the paraboloid (and in fact also for the sphere, by a slight modification of the argument); this was achieved in \cite{tao:parabola}.  We sketch the main idea of that paper here.

As in the last section, we can reduce matters to estimating a quantity such as
$$
\| \sum_{T_1 \in \T_1} \sum_{T_2 \in \T_2}  \psi_{T_1} \psi_{T_2} (1 - \chi_{B_{T_1} \cap B_{T_2}}) \|_{L^q(B(0,R))}$$
for some $1 < q < 2$.  It turns out in this case that one can obtain good bounds simply by interpolating between $L^1$ and $L^2$ bounds.  The $L^1$ bound is fairly trivial (using Cauchy-Schwarz to reduce to $L^2$ bounds on $\sum_{T_j \in \T_j} \psi_{T_j}$, which can be handled by orthogonality arguments), so we turn to problem of estimating the $L^2$ quantity:
$$
\| \sum_{T_1 \in \T_1} \sum_{T_2 \in \T_2}  \psi_{T_1} \psi_{T_2} (1 - \chi_{B_{T_1} \cap B_{T_2}}) \|_{L^2(B(0,R))}^2.$$
As is customary, we subdivide the large ball $B(0,R)$ into cubes $q$ of size $\sqrt{R}$.  The contribution of each cube $q$ is 
$$
\| \sum_{T_1 \in \T_1} \sum_{T_2 \in \T_2}  \psi_{T_1} \psi_{T_2} (1 - \chi_{B_{T_1} \cap B_{T_2}}) \|_{L^2(q)}^2.$$
Roughly speaking, we only need to consider pairs $T_1$, $T_2$ of tubes which pass through $q$ (because of the localization of $\psi_{T_1}$ and $\psi_{T_2}$, and such that $q$ is not contained in both $B_{T_1}$ and $B_{T_2}$.  For sake of argument, suppose that we only consider the terms where $q \not \subset B_{T_1}$.  Then we can rewrite the above expression as
$$
\| \sum_{T_1 \in \T'_1(q)} \psi_{T_1} \sum_{T_2 \in \T_2(q)} \psi_{T_2} \|_{L^2(q)}^2$$
where $\T_2(q)$ denotes all the tubes $T_2$ in $\T_2$ which intersect $q$, and $\T'_1(q)$ denotes all the tubes $T_1$ in $\T_1$ which intersect $q$ and for which $q \not \subset B_{T_1}$.  Note that the tubes in $\T'_1(q)$ must point in essentially different directions (since they all go through $q$, and are essentially distinct tubes), and similarly for $\T_2(q)$.

The function $\sum_{T_1 \in \T'_1(q)} \psi_{T_1}$ behaves roughly like the function $\frac{1}{R} (\chi_{\Omega'_1(q)} d\sigma_1)^\vee$, where $\Omega'_1(q)$ is the subset of the paraboloid whose unit normals lie within $1/R$ of the directions of one of the tubes in $\T'_1(q)$.  This can be seen by recalling the origin of these wave packets $\psi_{T_1}$, as Fourier transforms of functions on $S$ (or more precisely on $N_{1/R}(S)$; the discrepancy between the two explains the $\frac{1}{R}$ factor).  The function $\sum_{T_2 \in \T_2(q)} \psi_{T_2}$ is similarly comparable to the expression $\frac{1}{R} (\chi_{\Omega_2(q)} d\sigma_2)^\vee$ for a suitable set $\Omega_2(q)$.  Thus one is faced with an expression of the form
\begin{equation}\label{eq:22}
\| (\chi_{\Omega'_1(q)} d\sigma_1)^\vee (\chi_{\Omega_2(q)} d\sigma_2)^\vee \|_{L^2(\R^n)}^2,
\end{equation}
where we have discarded some powers of $R$ for sake of exposition, as well as the localization to $q$.  As observed previously, the restriction estimate $R^*(2 \times 2 \to 2)$ allows us to bound this quantity by something proportional to $|\Omega'_1(q)| |\Omega_2(q)|$.  Actually we can do a little better and refine this to, say, $|\Omega'_1(q)| |\Omega_2(q)|^2$; this comes from not discarding the constraint $\xi_2 \in \Omega_2$ in the argument immediately preceding \eqref{eq:crude}, and by exploiting the localization to $q$ more; we omit the details.  This is the type of bound used in Wolff's argument, in combination with the combinatorial arguments controlling the multiplicity of the tubes in $\T_1$ and $\T_2$ mentioned in the previous section, to obtain a sharp bilinear estimate in the case of the cone.  However, this bound is insufficient for the paraboloid case because of the failure of the tubes $T_1$ through a point to lie on a hypersurface.

Fortunately, this can be rectified by exploiting the gain inherent in \eqref{eq:crude}.  Indeed, by refusing to use \eqref{eq:crude} one can obtain a bound on \eqref{eq:22} which is proportional to
$$ |\Omega'_1(q)| |\Omega_2(q)| \sup_{\underline{\xi}_1, \underline{\xi}_4} |\Omega_2(q) \cap \pi(\underline{\xi}_1, \underline{\xi}_4)|.$$
This is similar to the bound of $|\Omega'_1(q)| |\Omega_2(q)|^2$ mentioned earlier, but is a little improved because one of the factors of $\Omega_2(q)$ is restricted to a hyperplane.  When one inserts this bound back into the coarse-scale combinatorial analysis of Wolff, this effectively allows us to restrict the tubes $T_2$ passing through a cube $q$ to be incident to a hyperplane.  This turns out to be a good substitute for the hypersurface localization property used in Wolff's argument, and is the key new ingredient which permits us to generalize the bilinear cone estimate to paraboloids (and by similar reasoning to other positively curved surfaces, such as the sphere).  

One interesting feature of this argument is that it introduces a non-trivial correlation between the fine-scale analysis and the coarse-scale analysis; one may speculate that future developments on these problems will deal with the fine-scale and coarse-scale aspects of the restriction operator in a more unified manner.

\end{document}